\documentclass[10pt]{article}
\usepackage[latin1]{inputenc}
\usepackage{amsmath,amssymb}
\usepackage{color}
\usepackage{multicol}
\usepackage{latexsym,epsfig}
\usepackage{amssymb,amsmath,amsfonts,amscd}
\usepackage{exscale}
\usepackage{ifthen}
\usepackage{longtable}
\usepackage{subfigure}
\usepackage{cite}
\usepackage{lscape}

\title{Iterated function systems with place dependent probabilities and  \\
application to the Diaconis-Friedman's  chain  on $[0,1]$}

\newtheorem{theo}{\textbf{Theorem}\ }
[section]
\newtheorem{lemma}[theo]{\textbf{Lemma}\ }

\newtheorem{defi}[theo]{\textbf{Definition}\ }

\begin{document}
\maketitle
\begin{center}
  Fetima   Ladjimi   $^($\footnote{Fetima Ladjimi,
Département de Mathématiques,  Universit\'e Mouloud MAMMERI
Tizi-Ouzou (UMMTO), 15000, Algérie -- email :  
ladjimi-fetima@ummto.dz }$^)$ \& Marc Peigné $^($\footnote{Marc
Peign\'e, LMPT, UMR 7350, Facult\'e des Sciences et Techniques,
Parc de Grandmont, 37200 Tours, France -- email :
peigne@univ-tours.fr}$^)$
$^($\footnote{ The author thanks the {\bf Vietnam Institute for Advanced Studies in Mathematics} (VIASM) in Ha Noi for its hospitality  during the first semester 2017}$^)$
 \vspace{2mm}
\end{center}

\begin{abstract}
   We study  Markov chains generated by  iterated Lipschitz functions systems
  with possibly place dependent probabilities. Under general conditions,
  we prove  uniqueness of the invariant probability measure for the associated Markov chain,
  by using  quasi-compact linear operators technics. We use  the same approach
   to describe the behavior of the Diaconis-Friedman's  chain  on $[0,1]$
   with possibly place dependent probabilities.

\end{abstract}



\section{Introduction}
We are interested in the  Markov  chain $(Z_n)_{n \geq 0}$
on $[0, 1]$
 introduced by P. Diaconis and D. Friedman in \cite{DF}. As it is described there,
  if the chain is at $ x$ at time $n$, it
selects at time $n+1$ one of the two intervals,  $[0, x]$ or $[x,
1]$ with equal probability $\frac{1}{2}$, and then moves  to a
random point $y$ in the chosen interval.

For $x \in ]0, 1[$, the transition
probability of the chain $(Z_n)_{n \geq 0}$ has a density $k(x,\cdot)$
with respect to the Lebesgue measure on $]0, 1[$ given by
$$\forall y \in ]0, 1[\qquad k(x,y)= \frac{1}{2}\times{1\over x}  {1\hspace*{-0.16cm}1}_{]0,x[}(y) + \frac{1}{2}
\times{1\over1-x}  {1\hspace*{-0.16cm}1}_{]x, 1[} (y).$$
Starting from $0$ (resp. $1$), the chain stays in $0$ (resp. $1$) with probability $\frac{1}{2}$ or moves with probability $\frac{1}{2}$ to a (uniformly chosen) random point in $]0, 1[$.

 It is shown in \cite{DF}  that it possesses a unique   invariant
probability measure $\nu$ on  $]0, 1[$;   this measure is the famous
``arcsine law'' which admits the density $f_{1\over 2}$
 with respect to the Lebesgue measure on $]0, 1[$  given by
 $$  \forall x \in ]0, 1[, \quad f_{1\over 2}(x)=\frac{1}{\pi \sqrt{x(1-x)}}{1\hspace*{-0.16cm}1}_{]0, 1[}(x).$$

 The same applies when the intervals $]0, x[$ and  $]x,
1[$ are chosen with the respective probabilities $p\in ]0, 1[$ and
$q=1-p$.  In this case, the invariant probability measure
 is the  Beta distribution $\mathcal B(q, p)$ of parameters $q$ and $p$
 whith  density  $ f_p$   defined by:
 $$\forall x \in ]0, 1[, \quad  f_p(x)= \frac{1}{\Gamma(p)\Gamma(q)}x
^{q-1}(1-x)^{p-1} {1\hspace*{-0.16cm}1}_{]0, 1[}(x).$$

 The transition  operator  $Q$ of  the chain $(Z_n)_{n \geq 0}$ is defined by:
 for any bounded
 Borel  function  $\varphi:[0, 1]\to \mathbb C$, 
 $$
 Q\varphi(0)= p\varphi(0)+q  \int_0^1 \varphi(y){\rm d}y, \qquad
 Q\varphi(1)= p  \int_0^1 \varphi(y){\rm d}y +q\varphi(1)
 $$
 and
$$
 \forall x \in ]0, 1[ \qquad Q\varphi(x) = {p\over x}\int_0^x \varphi(y){\rm d}y+ {q\over1-x}\int_x^1 \varphi(y){\rm d}y.
 $$
We may rewrite shortly $Q$ as follows: for any $x \in [0, 1]$,
\begin{equation}\label{noyau trans1}
Q\varphi(x)
= p \int_0^1 \varphi(tx) {\rm d} t+ q \int_0^1 \varphi(tx+1-t) {\rm d} t.
\end{equation}
This last expression shows that the
chain $ (Z_n) _ {n \geq 0} $ fits into the framework of iterated
random continuous functions. For any $t \in [0, 1]$,   let $H_t$ be the homothety $x \mapsto tx$ and  $A_t$ be the affine transformation
$x \mapsto tx+1-t$ and denote by $\mu$  the probability measure on  the space $ C([0,1], [0, 1])$ of continuous functions from $[0, 1]$ to  $[0, 1]$ defined by
$$\mu ( {\rm d}T)=p\int_0^1 \delta_{H_t}({\rm d}T) {\rm d}t+ q\int_0^1 \delta_{A_t}({\rm d}T) {\rm d}t,$$
where  $\delta_T$ is the Dirac masse at $T$. Equality
(\ref{noyau trans1})  may be restated as
$$Q\varphi(x)=\int_{ C([0,1], [0, 1])} \varphi(T(x)) \mu({\rm d}T).$$
Thus we may introduce a sequence  $(T_n)_{n \geq 1}$ of
independent random variables defined on a probability space
$(\Omega, \mathcal T, \mathbb P)$   with law   $\mu$ on $ C([0,1],
[0, 1])$. We have  $Z_n=T_n \cdots  T_1\cdot Z_0$; in other words,
the chain  $ (Z_n) _ {n \geq 0} $ is generated by iterating random
functions and its behavior is strongly connected to the
contraction  properties  of the  maps $H_t$ and $A_t, 0\leq t \leq
1$. We refer to $ (Z_n) _ {n \geq 0} $ as {\it the
Diaconis-Friedman's  chain.}

In  \cite{DF}, the authors  focus on the case when the  weights  $p$ and  $q$
 are depending on the  position $x$, and in particular when
   $p(x)=1-x$.  In the sequel, we  propose a  systematic examination
    of the general situation addressed by the two authors.

The study of Markov processes generated by  composition products
of random  independent functions  $ T_n $ has been the object of
numerous works for 50 years. When  the probabilities that govern the choice of
these transformations are spatially varying, the study of these
processes   escapes  the random walks framework. We refer
the reader to \cite{GR}, \cite {L}, \cite{MR} or \cite{ST} and
references there in, and  to   \cite{HH} or
\cite{P1}  for the approach via the theory of
quasi-compact operators. We use the terminology on
Markov chains as stated in \cite{R}.

\section{Iterated function systems}

     Let $(E, d)$ be  a metric compact space and denote $C(E, E)$ the space
of continuous functions   from  $E$  to $E$ endowed with the norm
$\vert\cdot\vert _\infty$ of the  uniform convergence on $E$. Let
$(T_n)_{n \geq 1}$ be a sequence of i.i.d random continuous
functions from $E$ to $E$ with distribution $\mu$.  The case when
the $T_n$ are   Lipschitz continuous from $E$ to $E$ is fruitful,
in particular   to use the so-called ``spectral gap property'',
based on the properties of contraction of the closed semi-group
$T_\mu$ generated by the support of $\mu$.

\subsection{ Iterated function systems with place independent probabilities}

 We denote
 $\displaystyle \mathbb L {\rm ip}(E, E) $ the space of
Lipschitz continuous functions from $E$ to $E$, i.e. of functions
$f: E\to E$ such that 
$$ [f]= \sup_{\stackrel{x, y
\in E}{x\neq y}} {d( f(x),f(y)) \over d(x,y)} < \infty,
$$
and we endow $\mathbb L {\rm ip}(E, E)$  with the norm $\Vert \cdot \Vert = \vert \cdot \vert_\infty+[\cdot]$.
  Let $(T_n)_{n \geq 1}$ be a sequence of
independent random functions  defined on a probability space
 $(\Omega, \mathcal T, \mathbb P)$,   with values in $\mathbb L {\rm ip}(E, E)$
 and common distribution $\mu$. We consider the Markov chain $(X_n)_{n \geq 0} $ on $E$,
 defined by:  for any $n\geq 0$,
 $$ X_{n+1}:= T_{n+1}(X_n)$$
 where $X_0$ is a fixed random variable  with values in $E$. One says that
  the chain $(X_n)_{n \geq 0} $ is
generated by   the {\it iterated function system}  $(T_n)_{n \geq
1}$. Its transition operator $P$  is defined by: for any bounded
Borel
 function  $\varphi: E \to \mathbb C$
 and   any  $x\in E$
$$ P\varphi(x)= \int_{\mathbb L {\rm ip}(E, E)} \varphi(T(x)) \mu({\rm d}T).$$
 The chain  $(X_n)_{n \geq 0} $ has the ``Feller property'',
i.e.  the operator $P$ acts on the space  $C(E)$ of continuous
functions from $E$ to $\mathbb C$. The maps $T_n$ being Lipschitz
continuous on $E$,  the operator $P$ acts also on the space of
Lipschitz continuous  from  $E$  to $\mathbb C$ and more generally
on the space  $ \mathcal   H_\alpha (E), 0<\alpha \leq 1$, of
$\alpha$-H\"older continuous functions from $E$ to $\mathbb C$,
defined by
$$
\mathcal   H_\alpha (E):= \{f\in C(E)\mid \Vert f\Vert_\alpha:=\vert f\vert_\infty +m_\alpha(f) <+\infty\}
$$
where $\displaystyle m_\alpha(f):=
\sup_{\stackrel{x, y \in E}{x\neq y}} {\vert f(x)-f(y)\vert \over
d(x,y)^\alpha} < \infty. $ 
Endowed, with the norm $\Vert \cdot
\Vert_\alpha$, the space $\mathcal   H_\alpha (E)$ is a Banach
space and the identity map from $C(E)$ to $\mathcal   H_\alpha
(E)$ is compact.

 The  behavior of the chain  $(X_n)_{n \geq 0} $ is closely related to
 the spectrum of $P$ on these spaces; under some ``contraction  in mean'' assumption
 on the  $T_n$, the restriction of the operator $P$ to  $\mathcal   H_\alpha (E)$ satisfies
 some spectral  gap property.
We first cite  the following theorem, due to Diaconis $\&$
Friedman  \cite{DF}; we  detail the proof  for the sake of
completeness.

\begin{theo}\label{muconstant}
Assume that there   exists  $\alpha \in ]0, 1]$
 such that
\begin{equation} \label{contraction}
r:= \sup_{\stackrel{x, y \in E}{x\neq y}}\int_{\mathbb L{\rm ip}(E, E)}\Bigl({ d(T(x),T(y)) \over d(x, y)}\Bigr)^\alpha \mu({\rm d}T) <1.
\end{equation}
Then there exists   on $E$ a unique $P$-invariant probability measure    $\nu$.  Furthermore, there exists
 constants $\kappa >0$ and $\rho \in ]0, 1[$    such that
\begin{equation} \label{vitesseExp}
\forall \varphi\in  \mathcal   H_\alpha (E),  \ \forall x \in E \quad  \vert P^n\varphi(x)-\nu(\varphi)\vert \leq \kappa \rho^n.
\end{equation}
\end{theo}
\noindent {\bf{Proof}}.
  The Feller operator $P$ is Markovian, thus  its spectral radius  $\rho_\infty(P)$
   in $C(E)$ equals $1$.
Furthermore,  $P$ acts on $\mathcal   H_\alpha(E) $ and for any function
$\varphi \in \mathcal   H_\alpha(E) $,  it holds
\begin{equation}\label{contracte}
 m_\alpha(P\varphi ) \leq
r\ m_\alpha(\varphi ),
\end{equation}  which yields
\begin{equation} \label{DFr}
\forall \varphi \in \mathcal   H_\alpha(E),\qquad
\Vert P\varphi \Vert_\alpha \leq r \Vert \varphi \Vert_\alpha
+\vert  \varphi \vert_\infty.
\end{equation}
Inequality  (\ref{DFr}) allows us to use the Ionescu-Tulcea and Marinescu theorem
for  quasi-compact operators. By Hennion's work \cite{H}, it implies that the essential
spectral radius of $P$ on $\mathcal   H_\alpha(E)$  is less than  $r$; in other words, any spectral
values  with modulus strictly larger then $r$ is an eigenvalue of $P$ with finite multiplicity
 and is  isolated in the spectrum of $P$.

  To prove the theorem, it is sufficient to control the peripheral spectrum of  $P$ on $\mathcal   H_\alpha(E)$. Let $\lambda $ an eigenvalue of  $P$ of modulus 1 and consider an  eigenfunction $f$ associated to $\lambda$.  For any $n \geq 1$, the equality   $  P^nf=\lambda^nf $   combined with (\ref{DFr})  yields
$$
m_\alpha(f)= m_\alpha(\lambda^n f)=m_\alpha(P^nf)\leq r^nm_\alpha(f)
$$
which implies $m_\alpha(f)=0$, since $0\leq r<1$.
Consequently, the  function  $f$ is  constant on $E$ and  $\lambda = 1$.
Thus, the operator $P$  on  $\mathcal   H_\alpha(E)$  can be decomposed  as,
\begin{equation} \label{Pi+R}
P= \Pi +R
\end{equation}
where

(i) the operator $\Pi  $ is  the projector from $\mathcal   H_\alpha(E)$ to the eigenspace  $\mathbb C\cdot {1\hspace*{-0.16cm}1}$ associated to the eigenvalue $1$,

(ii) $R$ is an operator with spectral radius $\rho $ for some $\rho \in [0, 1[$,

(iii)  $\Pi R= R \Pi = 0$.

In particular, for any $\varphi\in \mathcal   H_\alpha(E)$, the sequence $(P^n\varphi)_{n \geq 0}$ converges to $\Pi(\varphi){1\hspace*{-0.16cm}1}$; thus, there exists  on $E$ a unique invariant probability measure  $\nu$ and the projector   $\Pi$ may be written as
$\Pi: \varphi \mapsto \nu(\varphi) {1\hspace*{-0.16cm}1}$.  Inequality
(\ref{vitesseExp})  follows from decomposition (\ref{Pi+R}).

\rightline{\rightline{$\Box$}}

 \noindent {\bf Application to the
Diaconis-Friedman's chain for $p$ fixed in $]0, 1[.$}

Inequality (\ref{contraction})  holds  with $r= {1\over 1+\alpha}
$ since  in this case $m(H_t)=m(A_t)=t$ for  any $0\leq t\leq 1$.
Hence
\begin{eqnarray*}
\sup_{\stackrel{x, y \in [0, 1]}{x\neq y}}\int_{\mathbb L{\rm
ip}([0, 1], [0, 1])}\Bigl({ d(T(x),T(y)) \over d(x,
y)}\Bigr)^\alpha \mu({\rm d}T)  &\leq& p\int_0^1  m(H_t) ^\alpha
{\rm d}t+q\int_0^1  m(A_t) ^\alpha {\rm d}t
\\
&
=&
 \int_0^1 t^\alpha {\rm d}t={1\over 1+\alpha}.
\end{eqnarray*} 

Thus, the chain  $(Z_n)_{n \geq 0}$ admits  an unique  invariant probability  measure   on $[0,1]$, this measure being the Beta distribution $\mathcal  B(q, p)$.

\rightline{\rightline{$\Box$}}

\subsection{Iterated function systems with  spacial dependant increments probabilities}

It this section, we replace the measure $\mu$ by a collection
$(\mu_x)_{x \in E}$ of probability measures on $E$, depending
continuously on $x$. We consider the Markov chain $(X_n)_{n \geq
0}$ on $E$ whose  transition kernel $P$ is given by: for any
bounded  Borel function  $\varphi: E \to \mathbb C$
 and  any  $x\in E$,
$$
P\varphi(x)= \int_{\mathbb L {\rm ip}(E, E)} \varphi(T(x)) \mu_x({\rm d}T).
$$
 First, we introduce  the following definition.

\begin{defi}
  A sequence  $(\xi_n)_{n \geq 0}$ of  continuous functions from $E$ to $E$ is  a contracting sequence  if there exist  $x_0 \in E$ such that
  $$
\forall x \in E \quad \lim_{n \to +\infty} \xi_n(x)= x_0.$$
\end{defi}
The following statement is a generalization of Theorem \ref{muconstant}.
\begin{theo}\label{muvarie}
Assume that there   exists $\alpha \in (0,1]$ such that
\begin{enumerate}
\item [{\bf H1.}]  $r:=\displaystyle \sup_{\stackrel{x, y \in E}{x\neq y}}\int_{\mathbb L{\rm ip}(E, E)}\Bigl({ d(T(x), T(y))\over d(x, y)}\Bigr)^\alpha\mu_x({\rm d}T) <1,
 $
\item [{\bf H2.}]   $R_\alpha:= \displaystyle
\sup_{\stackrel{x, y \in E}{x\neq y}}{\vert \mu_x-\mu_y\vert \over d (x, y)^\alpha}<+\infty,
$
\item [
{\bf H3.}] there exists  $\delta >0$ and a  probability measure
$\mu$ on  $E$  such that
\begin{equation} \label{minoration}
\forall x \in E\qquad \mu_x\geq \delta \mu
\end{equation}
and the closed  semi-group   $T_\mu$  generated by the support
  $S_\mu$  of $ \mu$ possesses a contracting sequence.
\end{enumerate}
Then, there exists   on $E$ a unique $P$-invariant probability measure    $\nu$; furthermore, for some constants $\kappa >0$ and $\rho \in ]0, 1[$, it holds
\begin{equation}\label{vitesseExpbis}
\forall \varphi \in \mathcal   H_\alpha(E), \ \forall x \in E \quad \vert P^n\varphi(x)-\nu(\varphi)\vert \leq \kappa\rho^n.
\end{equation}
\end{theo}
 \noindent {\bf Remark.} Hypothesis {\bf H1} means that the maps
$T$ satisfy some contraction property ``in mean'', with respect to
each measure $\mu_x$. Nevertheless, the measures $\mu_x$  may be
singular versus another; this implies that, starting from two
different points, the maps which govern the transition may be
totally different and it becomes quite impossible to control their
 common evolution. Thus, hypothesis {\bf H3} is
useful to fill up this gap.

\noindent {\bf Proof.}  The  operator $P$ acts on $C(E)$, with
spectral radius $1$ since it is Markovian. It also acts on
$\displaystyle \mathcal   H_\alpha(E) $; indeed, for any function
$\varphi \in    \mathcal   H_\alpha(E)  $ and any $x, y \in E$, it
holds
\begin{eqnarray*}
\vert P\varphi(x)-P\varphi(y)\vert &\leq&  \int_{\mathbb L{\rm ip}(E, E)}\vert \varphi(T(x))-\varphi(T(y))\vert   \mu_x({\rm d}T) +\vert\varphi\vert_\infty  \int_{\mathbb L{\rm ip}(E, E)}\vert\mu_x-\mu_y\vert ({\rm d}T).
\end{eqnarray*}
Hence
\begin{equation}\label{malpha}
m_\alpha(P\varphi)\leq r m_\alpha(\varphi)   + R_\alpha\vert \varphi\vert _\infty
\end{equation}
which readily yields  
\begin{equation}\label{DF}\Vert
P\varphi\Vert_{\alpha}\leq r\Vert \varphi\Vert_\alpha  +
(1+R_\alpha)\vert \varphi\vert _\infty.
\end{equation} 
 Thus, by \cite{H}, the operator $P$ is
quasi-compact  on $   \mathcal   H_\alpha(E)$; its spectral radius on
 $   \mathcal   H_\alpha(E)$ equals the modulus of a dominant eigenvalue, thus  is less than the one of $P$ on $C(E)$, that is $1$. To control the
 peripheral spectrum, the argument differs then from
the one used to prove Theorem \ref{muconstant}:  property
(\ref{contracte}) does not hold  here and  inequality
(\ref{malpha}) is   much weaker. We  get use of the two following
lemmas, valid under  hypotheses {\bf  H1}, {\bf  H2}
 and {\bf  H3}.

\begin{lemma} \label{harmonic} Let $h \in \mathcal   H_\alpha(E) $ such that $Ph=h$.
For any $x \in E$, the sequence  $(h(X_n))_{n \geq 0}$ is a
bounded martingale  on the space $(\Omega, \mathcal F, \mathbb
P_x)$,  where $\mathbb P_x$ denotes the conditional probability
$\mathbb P(\cdot / X_0=x)$. It  converges $\mathbb P_x$-a.s. and
 in $\mathbb L^1(\Omega, \mathbb P_x)$  to  a random variable  $H^\infty$ and it holds
\begin{equation}\label{martingale}
\forall n \geq 0 \quad h(x)= \mathbb E_x(h(X_n))=  \mathbb
E_x(H^\infty).
\end{equation}
Furthermore, for any $\xi \in T_\mu$,
\begin{equation}\label{convergence}
H^\infty= \lim_{n \to +\infty}  h(\xi\cdot X_n) \quad \mathbb P_x-{\rm a.s.}
\end{equation}

\end{lemma}
 Proof of Lemma \ref{harmonic}. The
function $h$ is $P$-harmonic and bounded;  the first  assertion
and equality (\ref{martingale}) follow. Let us now prove
(\ref{convergence}). First, let us fix    positive integers $n$
an $q$ and set
$$u_{n, q}(x):= \mathbb E_x\Bigl(\bigl\vert h(X_{n+q})-h(X_{n})\big\vert^2\Bigr).$$
  From the martingale equality, for any  $N\geq 1$, it holds
$$ \sum_{n=1}^N u_{n, q}(x)= \sum_{n=1}^N
\mathbb E_x\Bigl(\big\vert h(X_{n+q})\big\vert^2\Bigr)-
\sum_{n=1}^N\mathbb E_x\Bigl(\big\vert h(X_{n})\big\vert ^2\Bigr) \leq 2q
\big\vert h\big\vert_\infty ^2. $$
Hence, $\displaystyle \sum_{n =1}^{+\infty} u_{n, q}(x) <+\infty$ and
\begin{equation}  \label{sumvaleurpropre}
 \sum_{n =1}^{+\infty}
\mathbb E_x\Bigl(  \int_{\mathbb L{\rm ip}(E, E)^q}\big\vert  h(T_{q}\cdots T_1\cdot X_n)-h( X_n)\big\vert ^2 \mu_{X_n}({\rm
d}T_1)\cdots \mu_{T_{q-1}\cdots T_1\cdot X_n}({\rm d}T_q)\Bigr)
 <+\infty.
\end{equation}
Consequently, using    H3,
$$
\sum_{n =1}^{+\infty}
\mathbb E_x\Bigl(  \int_{\mathbb L{\rm ip}(E, E)^q}\big\vert  h(T_{q}\cdots T_1\cdot X_n)-h( X_n)\big\vert ^2 \mu({\rm
d}T_1)\cdots \mu({\rm d}T_q)\Bigr)
 <+\infty $$
and
$$
\int_{\mathbb L{\rm ip}(E, E)^q}\mathbb E_x \Bigl( \sum_{n =1}^{+\infty}\big\vert
 h(T_{q}\cdots T_1\cdot X_n)-h( X_n)\big\vert^2 \Bigr) \mu ({\rm
d}T_1)\cdots \mu ({\rm d}T_q)
<+\infty.$$
 For any $q \geq 1$   and  $\mu^{\otimes q}$-almost all $T_1, T_2, \cdots, T_{q}$, the sequence
 $\displaystyle \left( h(T_{q}\cdots T_1 \cdot X_n)
-h(X_n)\right)_{n \geq 1}$ converges $\mathbb P_x$-a.s.  to $0$.
We conclude by  a density argument.

\rightline{$\Box$}

 Similarly, one may prove the following lemma, which is of
interest  to control   the other modulus 1  eigenvalues of $P$ in
$ \mathcal   H_\alpha(E) $.
\begin{lemma} \label{valeurpropre} Let $\phi \in \mathcal   H_\alpha(E) $ such that $P\phi= \lambda \phi$ where $\lambda$ is a complex number of   modulus  1. For any $x \in E$, the sequence  $(\lambda^{-n}\phi (X_n))_{n \geq 0}$ is a bounded  martingale; it  converges $\mathbb P_x$-a.s and in  $\mathbb L^1(\Omega, \mathbb P_x)$ to a random variable $\Phi^{\infty}$ and we have
\begin{equation}\label{martingaleET}
\forall n \geq 0 \quad \phi(x)= \mathbb E_x(\lambda^{-n}\phi(X_n))=  \mathbb E_x( \Phi^{\infty})
\end{equation}
 Furthermore, for any  $q\geq 1$ and any transformations   $T_{1}, \cdots, T_{q}$ on the support $S_\mu$ of $\mu$,  one has
\begin{equation}\label{semigroupeET}
\Phi^\infty= \lim_{n \to +\infty}  \lambda^{-(n+q)}\phi(T_q\cdots T_1\cdot X_n) \quad \mathbb P_x-{\rm a.s.}
 \end{equation}
\end{lemma}

Let us first  prove that the  $P$-harmonic functions in $ \mathcal   H_\alpha(E)
$ are constant. Let  $h \in \mathcal   H_\alpha(E) $ such that
 $Ph=h$. According to   Lemma \ref{harmonic},
for any $x \in E$,  there exists a set   $\Omega_x\subset \Omega $  of full measure with respect
to   $\mathbb P_x $ such that,  for any $ \omega \in \Omega_x$
and  any transformation $\xi \in T_\mu$, the sequences
 $(h(X_n(\omega)))_{n\geq 0}$ and $(h(\xi\cdot X_n(\omega)))_{n\geq 0}$
 converge  to $H^\infty(\omega).$

Let  $(\xi_k)_{k\geq 0}$ be a contracting sequence  in $T_\mu$,
with   limit point $x_0\in E$.  Since  $h  $ is continuous on $E$,
for any $\omega \in \Omega_x$,  any cluster value $ x_\omega$ of
$(X_n(\omega))_{n \geq 0}$ and  any $k \geq 0$,
 $$H^\infty(\omega)= h(x_\omega)= h(\xi_k(
x_\omega)).$$ 
 Letting $k\rightarrow +\infty$, it yields
 $ H^\infty(\omega)= h(x_0)$ and    thus  $h(x)= h(x_0)$, by   (\ref{martingale}). Finally, the
 bounded $P$-harmonic functions in $   \mathcal   H_\alpha(E) $ are constant.

 Using Lemma \ref{valeurpropre}, we prove that the peripheral
spectrum of $P$ is reduced to $1$. Let $(n_l)_{l\geq 0}$ be a
fixed  sequence of integers  such that $\displaystyle
\lim_{l\to+\infty}\lambda^{-n_l}=1$ and $(\xi_k)_{k\geq 0} $ be a
contracting sequence on $T_\mu$, with limit point  $x_0$. \\ For
all integer $q\geq 1$,  the set  $\{T= T_q\cdots T_1 \ \mid\ T_1
\ldots, T_q\in S_\mu\}$ is  dense in $T_\mu$. Without loss of
generality, we assume that  any function $\xi_k$ can be decomposed
as a  product   $T_{q_k} \cdots T_1$, with $T_i \in S_\mu, 1\leq
i\leq q_k$. By  Lemma \ref{valeurpropre},  there exists
 $\Omega_x\subset \Omega$, $\mathbb P_x(\Omega_x)=1$, such that,
  for any $ \omega \in \Omega_x$ and $ k\geq 0$,
  the sequences  $\Bigl(\lambda^{-n_l}\phi(X_{n_l} (\omega)\Bigr)_{l\geq 0}$ and
$\Bigl(\lambda^{-(n_l+q_{k})}\phi(\xi_k\cdot X_{n_l}
(\omega))\Bigr))_{l\geq 0}$ converge to the same limit
$\Phi^\infty(\omega)$.   Let us choose sequences of
 integers  $(\varphi(l))_{l \geq 0} $  (depending  on $\omega$) and $(\psi(k))_{k \geq 0} $
  (which does not depend on $\omega$) such
  that $(\lambda^{-n_{\varphi(l)}}X_{n_{\varphi(l)}}(\omega))_{l \geq 0}$
   and $(\lambda^{-q_{\psi(k)}})_{k \geq 0}$ converge  resp. to
$x_\omega \in E$  and $e^{i\beta}, \beta \in \mathbb R$. Equalities (\ref{martingaleET}) and  (\ref{semigroupeET})  yield
$$\Phi^\infty(\omega)= \phi(x_\omega)= e^{i\beta} \phi(x_0)\quad {\rm and} \quad
\phi(x)=\mathbb E_x(\Phi^\infty)=e^{i\beta} \phi(x_0).$$
Eventually, the function   $\phi$ is constant on  $E$  and
$\lambda =1$.

\rightline{$\Box$}
\section{The  Diaconis-Friedman's chain}

 This section deals with the  Diaconis-Friedman's chain $(Z_n)_{n
\geq 0}$   on $E=[0, 1]$ described in  the introduction; we assume
that  the weights $p$ and $q$    vary with  $x\in [0, 1]$. The
transition  operator $Q$ of $(Z_n)_{n \geq 0}$ is given by:
 for any bounded  Borel function $\varphi:[0, 1] \to
\mathbb C$,
$$ Q\varphi(x)= p(x) \int_0^1 \varphi(tx) {\rm d} t+ q(x) \int_0^1 \varphi(tx+1-t) {\rm d} t.
$$

For  $x \in [0, 1]$, let $\mu_x$ be the probability  measure on
the space $\mathbb L {\rm ip}([0, 1], [0, 1])$   of Lipschitz continuous functions
from $[0, 1]$ into $[0, 1]$, defined by
\begin{equation}\label{diaconismux}
\mu_x ( {\rm d}T)=p(x)\int_0^1 \delta_{H_t}({\rm d}T) {\rm d}t+ q(x)\int_0^1 \delta_{A_t}({\rm d}T) {\rm d}t.
\end{equation}
 Then the  transition operator $Q$  may be rewritten as
 $$Q\varphi(x)=\int_{\mathbb L {\rm ip}([0, 1],[0,
1])} \varphi(T(x)) \mu_x({\rm d}T)= p(x)\int_0^1 \varphi(tx) {\rm
d}t+ q(x)\int_0^1 \varphi(tx+1-t) {\rm d}t.$$
   Let us first
consider explicit examples.

\begin{enumerate}
\item When $p(x) = x$,  the chain $(Z_n)_{n \geq 0}$ is a
sequence of independent random  variables  of  uniform
distribution on $[0, 1]$; thus,  its unique invariant measure is
the uniform distribution on $[0, 1]$.

\item When $p(x) = 1-x$,  the points   $0$ and $1$ are  absorbing points for $(Z_n)_{n \geq 0}$.
 Hence, the Dirac measures at $0$ and $1$   are  $Q$-invariant. The following theorem states that
  these two measures  are the only   ergodic probability measures on $[0, 1]$ and
   that  $(Z_n)_{n\geq 0}$  converges  $\mathbb P_x$-a.s.
    to a random variable $Z_\infty$ with values in $\{0,1\}$.

\item Assume that  $p\in{ \mathcal H}_\alpha[0, 1] $ satisfies
\begin{equation}\label{pstrictementpositif}
 \forall x \in [0, 1], \quad  p(x)>0
 \end{equation}
(or in a symmetric way, $ p(x)<1$ for any  $x \in [0, 1]$).
In this case, the chain $(Z_n)_{n \geq 0}$ admits a unique $Q$-invariant probability measure on
  $[0, 1]$.  This is  a direct consequence of Theorem   \ref{muvarie}. Indeed,
   hypotheses {\bf H1}, {\bf H2} and {\bf H3}  of Theorem \ref{muvarie} hold:
\begin{enumerate}
\item Hypothesis {\bf H1}.  For any $x, y \in [0, 1], x\neq y$,
$$\sup_{\stackrel{x, y \in [0, 1]}{x\neq y}}\int_{\mathbb L{\rm ip}([0, 1], [0, 1])}\left({ \vert
T(x)-T(y)\vert \over \vert x-y\vert} \right)^\alpha \mu_x({\rm
d}T) \leq    {1\over 1+\alpha}.
$$

\item  Hypothesis   {\bf H2}. For any $x, y \in [0, 1], x\neq y$,
$${\vert \mu_x-\mu_y\vert \over \vert x-y\vert^\alpha}\leq {\vert p(x)-p(y)\vert \over \vert x-y\vert^\alpha}\Big\vert \int_0^1 \delta_{H_t}{\rm d}t\Big\vert
+{\vert q(x)-q(y)\vert \over \vert x-y\vert^\alpha}\Big\vert \int_0^1 \delta_{A_t}{\rm d}t\Big\vert \leq 2m_\alpha(p).
$$
\item  Hypothesis  {\bf H3}.  For any $x \in [0, 1]$ it holds $ \mu_x \geq \delta \mu $
with $\displaystyle \delta := \inf_{x \in [0, 1]} p(x)>0$ and $\displaystyle \mu= \int_0^1 \delta_{H_t}{\rm d}t$.
 The constant function $H_0: x \mapsto 0$ belongs to the support of   $\mu$; hence,
 the semi-group $T_\mu$ contains a contracting sequence, with limit point $0$.
\end{enumerate}
 If   $p(0)=1$, the Dirac mass at $0$ is the unique invariant
probability  measure   for $(Z_n)_{n \geq 0}$. When $p(0)<1$, one
can prove that  the unique invariant probability measure  for
$(Z_n)_{n \geq 0}$ is absolutely continuous with respect to the
Lebesgue measure (see Theorem \ref{GENERALDF} below).
\end{enumerate}
  If $p$ and
 $ q $ are both strictly positive on $[0, 1]$,  by using the approach
 developed in \cite{DF}, we may  prove  that the unique invariant probability measure
 for $ (Z_n) _ {n \geq 0} $ is absolutely continuous with respect to the
 Lebesgue measure. This property holds  as soon as    $ p (0) <1 $ and $ q (1) <1 $.
   Let us emphasize that  the strict positivity of $p$ or  $q$ is sufficient
    to ensure the unicity of an invariant probability measure but it is  a too
    strong condition. These remarks lead  to  the following statement, which is
    not a  direct consequence of  Theorem \ref{muvarie} but whose  proof is strongly inspired.

 \begin{theo}\label{GENERALDF}
 Let  $(Z_n)_{n \geq 0}$ be the Diaconis-Friedman's chain on $[0,
1]$ with weight  functions $ p $ and $ q $ in $\mathcal
H_\alpha[0,1]$. Then, one of the 3 following options holds.
\begin{enumerate}
\item If $p(0)<1$ and $q(1)<1$,  then  there exists on $[0, 1]$ an unique $Q$-invariant
 probability measure $\nu_p$. Furthermore, this measure is absolutely continuous
 with respect  to the  Lebesgue measure on $[0, 1]$ with density $f_p$ given by:
$$
\forall x \in [0, 1]\quad f_p(x)= C \exp\Bigl(\int_x^{1\over 2}{p(y)\over y} {\rm d}y+ \int_{1\over 2}^{x}{q(y)\over 1-y} {\rm d}y\Bigr)
$$
where $C$ is a normalization constant. At last, there exist
constants $\kappa >0$ and $\rho \in [0, 1[$ such that
$$
\forall \varphi \in \mathcal   H_\alpha[0,1], \ \forall x \in [0, 1]\quad \vert Q^n\varphi(x)-\nu_p(\varphi)\vert \leq \kappa \rho^n\Vert \varphi\Vert_\alpha.
$$
\item If $p(0)=1$ and $q(1)<1$, then  the  Dirac measure
 $\delta_0$ is the unique $Q$-invariant probability measure
  on $[0, 1]$. Furthermore, there exist  constants
  $\kappa >0$ and $\rho\in [0, 1[$ such that
   $$\forall \varphi \in \mathcal   H_\alpha[0,1], \ \forall x \in [0, 1] \quad \vert Q^n\varphi(x)-\varphi(0)\vert \leq  \kappa \rho^n \Vert \varphi \Vert_\alpha.$$
(A similar statement holds when $p(0)<1$ and  $q(1)=1$).
\item   If $p(0)=1$ and $q(1)=1$,
then the invariant  probability measures of $(Z_n)_{n \geq 0}$
  are the convex combinations
of  $\delta_0$ and  $\delta_1$. Furthermore, for any $x \in [0, 1]$,
 the chain $(Z_n)_{n \geq 0}$ converges
 $\mathbb P_x$-a.s. to a random   variable
  $Z_\infty$ with values in $\{0, 1\}$; the   law  of $Z_\infty$ is given by
   $$ \mathbb P_x(Z_\infty=0)=1-h(x)\quad  and  \quad \mathbb P_x(Z_\infty=1)= h(x), $$
    where  $h$ is  the  unique function  in $\mathcal   H_\alpha[0, 1]$ such that  $Qh=h$ and
      $h(0)=0, \ h(1)=1$.
At last, there exist $\kappa >0$ and $\rho\in [0, 1[$ such that
$$ \forall \varphi \in \mathcal   H_\alpha[0,1], \ \forall x \in [0, 1] \quad \vert Q^n\varphi(x)-(1-h(x))\varphi(0)-h(x)\varphi(1)\vert \leq  \kappa \rho^n \Vert \varphi\Vert_\alpha.$$
\end{enumerate}

\end{theo}

\noindent {Proof}. First,   let us consider the adjoint  operator  $Q^*$  of $Q$ in
$\mathbb L^2[0, 1]$, defined by: for any $ \varphi, \psi  \in \mathbb L^2[0, 1]$,
$$\int_0^1 \varphi(x) Q\psi(x) {\rm d}x= \int_0^1 Q^*\varphi(x)  \psi(x) {\rm d}x.$$
A straightforward computation yields to the following expression:
\begin{equation}\label{adjoint}
\forall \varphi \in \mathbb L^2[0, 1], \forall x \in [0, 1] \qquad  Q^*\varphi(x):= \int_0^x{q(t)\over 1-t}\varphi(t) {\rm d}t+
 \int_x^1{p(t)\over  t}\varphi(t) {\rm d}t.
\end{equation}
Notice that   (\ref{adjoint}) is valid  for any  Borel function $\varphi  \in \mathbb L^1[0, 1]$. Furthermore, if $\varphi  \in \mathbb L^1[0, 1]$ is non negative and satisfies the equality   $Q^*\varphi= \varphi$, then  the measure with density $\varphi$ with respect to the Lebesgue measure on $[0, 1]$  is  $Q$-invariant.

Assume for a while that  $\varphi$ is  differentiable on $]0, 1[$; the equation $Q^*\varphi= \varphi$  yields
$$\forall x\in ]0, 1[ \quad \varphi'(x)= \Bigl({q(x)\over 1-x} -{p(x)\over x}  \Bigr)\varphi(x), $$
hence $\displaystyle \varphi(x)=   \exp\Bigl(\int_x^{1/2}
{p(t)\over t} {\rm d}t+ \int_{1/2}^x {q(t)\over1-t} {\rm d}t
\Bigr) $ up to a multiplicative constant. This function $\varphi$
is integrable with respect to the Lebesgue measure on $[0, 1]$ if
and only if $p(0)<1$ and $q(1)<1$;  in this case, we set
\begin{equation}\label{densite}
f_p: x \mapsto  \displaystyle   {1\over C_p} \exp\Bigl(\int_x^{1/2}
{p(t)\over t} {\rm d}t+ \int_{1/2}^x {q(t)\over1-t} {\rm d}t
 \Bigr)
 \end{equation}
 with $C_p:= \displaystyle \int_0^1\exp\Bigl(\int_x^{1/2}
{p(t)\over t} {\rm d}t+ \int_{1/2}^x {q(t)\over1-t} {\rm d}t
 \Bigr){\rm d}x$. The probability measure $\nu_p$ on [0, 1] with density $f_p$
  with respect to the  Lebesgue measure  on $[0, 1]$  is
   $Q$-invariant.

\noindent  Now, we come back to the proof of  Theorem
\ref{GENERALDF} and decompose the argument into 3 steps.

\vspace{2mm}

\noindent  {\it  Step1-  Quasi-compacity of the operator $Q$ on
$\mathcal   H_\alpha[0,1]$}
\vspace{2mm}

The operator $Q$ is  non negative, bounded on $  \mathcal   H_\alpha [0, 1]$ with  spectral radius
 $1$. Furthermore,
$$\forall \varphi \in \mathcal   H_\alpha[0,1], \qquad \Vert Q\varphi \Vert_\alpha \leq {1\over \alpha+1} \Vert \varphi\Vert_\alpha +(1+2m_\alpha(p)) \vert \varphi\vert_\infty.$$
Hence, by \cite{H}, the operator  $Q$ is quasi-compact on $  \mathcal   H_\alpha [0, 1]$.

\vspace{2mm}

\noindent  {\it   Step  2- Description of the characteristic space of $Q$ corresponding to   $\lambda = 1$}

\vspace{2mm}  We use here a general result of \cite{Her}, based on
the notion  of {\it  absorbing compact set.} A compact subset $ K
$ of $ [0, 1] $ is said to be {\it$Q$-absorbing}  if
$Q1_{[0,1]\setminus K}(x)=0$ for any $x\in K$. It is {\it minimal}
when it does not contain any proper absorbing compact subset. The
condition $p(x)>0$ ensures that
\begin{equation}\label{*}
Q(x, I)>0 \quad {\rm for\  any\  closed \ interval} \  I\subset [0, x] {\rm  \ not \ reduced\  to \ a \ single \ point.}
\end{equation}
Similarly, the  condition $q(x)>0$ implies
\begin{equation}\label{**}
Q(x, I)>0 \quad {\rm for\  any\  closed \ interval} \  I\subset [x,1] {\rm  \ not \ reduced\  to \ a \ single \ point.}
\end{equation}

There are four cases to explore.
\begin{enumerate}
\item \  \underline{\it  $ q(0)>0$ and $ p(1)>0$}

 In this case, the interval $[0, 1]$ is the unique (and thus
minimal)  $Q$-absorbing compact set.  To prove this, we fix a
compact and proper subset $ K$ of $[0, 1] $ ; we have to find a
point   $x_0 \in K$ such that $Q(x_0, [0, 1]\setminus K)>0$. There
are 3 sub-cases to consider.
\begin{enumerate}
\item[(a)]  \underline{$0 \notin K$}

 Assume that
 $q(x)=1$  for any  $x \in
K$. The condition  $p(1)>0$ implies $q(1)<1$, so that $1\notin K$; thus, there
exist  $\epsilon>0$
 such that $K\subset [0, 1-\epsilon]$. Consequently,
for any $x  \in K$,
  $$Q(x, [0,1]\setminus K) \geq Q(x, ] 1-\epsilon, 1]) >0$$
  which means that $ K$ is not absorbing. Contradiction.

Consequently, there exists  $x_0  \in K$ such that $p(x_0 )>0$; if $\epsilon>0$ is  such that $K\subset [\epsilon, 1]$, then,
$$Q(x_0, [0,1]\setminus K) \geq Q(x_0,[0, \epsilon[) >0.$$

\item[(b)] \ \underline{$1 \notin K$}

The same argument holds,   exchanging the role of
$0$ and $1$.

\item[(c)] \ \underline{$0\in K $ and $1\in K$}

In this case, we can set  $x_0=0$. Indeed, let us fix  $x'\in ]0, 1[\setminus K$ and $\epsilon'>0$  such that
$]x'-\epsilon', x'+\epsilon'[ \subset [0, 1]\setminus K  $ and notice that
$$Q(0, [0, 1]\setminus K)\geq Q(0, ]x'-\epsilon', x'+\epsilon'[)>0.
$$
\end{enumerate}

 \item  \underline{ $ q(0)=0$ and $  p(1)>0$}
 In this case, the set  $\{0\}$ is invariant and  is the unique
$Q$-absorbing  minimal  compact  set. Indeed,   there exists $x
\in K$ such that $p(x) >0$. Otherwise, the function $q$ equals $1$
on $K$; by (\ref{**}), it follows that $[y, 1] \subset K$ for any
$y \in K$. Consequently  $1\in K$ and $q(1)=1$, which contradicts
the condition $  p(1)>0$. Applying  (\ref{*}), it yields
  $[0, x]\subset K$ and in particular $\{0\}\subset K$.

\item  \underline{     $ q(0)>0$ and $  p(1)=0$}

In this case, the
unique $Q$-absorbing  minimal  compact  set is $\{1\}$. The proof is similar to the previous  case, exchanging the role of
$0$ and $1$.

\item \underline{  $ q(0)=0$ and $  p(1)=0$}

The sets  $\{0\}$ and $\{1\}$ are the only minimal   absorbing
compact sets. \end{enumerate}  We apply Theorem 2.2 in \cite{Her}
to  conclude  that   the  eigenvalue $1$  has index $1$ in $
\mathcal   H_\alpha [0, 1]$: in other words,  the characteristic
subspace   of $Q$ associated to  $1$  equals ${\rm Ker}(Q-Id)$.
Therefore, we may  apply Theorem 2.3 in \cite{Her} to each of the
four cases explored above.
\begin{enumerate}
\item If $ q(0)>0$ and $ p(1)>0$,  then ${\rm Ker}(Q-Id)=\mathbb C\cdot {1\hspace*{-0.16cm}1}$; in this case, the unique $Q$-invariant probability measure   on $ [0, 1] $ is absolutely continuous with respect to the Lebesgue measure on $[0, 1]$, with density   $ f_p$.
\item  If $ q(0)=0$ and $ p(1)>0$,  then ${\rm Ker}(Q-Id)=\mathbb C\cdot {1\hspace*{-0.16cm}1}$   and the  Dirac mass $\delta_0$ is the unique $Q$-invariant probability measure   on $ [0, 1]$.
\item If $ q(0)>0$ and $ p(1)=0$, then    ${\rm Ker}(Q-Id)=\mathbb C\cdot {1\hspace*{-0.16cm}1}$ and the  Dirac mass $\delta_1$ is the unique $Q$-invariant probability measure   on $ [0, 1].$
\item If $ q(0)=0$ and $ p(1)=0$, there exists a positive harmonic function $h$ such that  $h (0)=0$ and $h (1)=1$;  the  space ${\rm Ker}(Q-Id)$  has  dimension $2$ and  equals  $\mathbb C\cdot {1\hspace*{-0.16cm}1}\oplus \mathbb C \cdot h$. The $Q$-invariant probability measure   on $ [0, 1]$ are the convex combinations of $\delta_0$ and $ \delta_1$.
\end{enumerate}

\vspace{2mm}

\noindent \underline{\it  Step 3- Control of the peripheral spectrum  of $Q$ in $  \mathcal   H_\alpha [0, 1]$}\ \\

\vspace{2mm}
  We use here Lemma \ref {valeurpropre} and apply    the same technics as in the previous discussion.

  Let $\lambda \in \mathbb C$ with modulus 1 and $\phi \in \mathcal   H_\alpha[0,1] $
  such that $Q\phi= \lambda \phi$. For any  $x \in [0, 1]$, the sequence
   $(\lambda^{-n}\phi(X_n))_{n \geq 0}$  is a bounded  martingale   in

   $(\Omega, \mathcal F, \mathbb P_x)$, thus it converges  $\mathbb P_x$-a.s to  a bounded
    random variable $\Phi^\infty.$
 We use inequality (\ref{sumvaleurpropre}) first with $q=1$ and
then $q=2$;    there exists   $\Omega_x \subset \Omega, \mathbb
P_x(\Omega_x) =1,$  and  $I_0\subset [0, 1]$  of Lebesgue measure
$1$ such that, for any $\omega \in \Omega_x$ and any  $s, t \in
I_0$, it holds
\begin{equation}\label{AA}
\lim_{n \to +\infty} \big\vert   \phi(Z_n(\omega))-\lambda^{-1}\phi(H_s\cdot Z_n(\omega))\big\vert^2 p(Z_n(\omega))=0,
\end{equation}
\begin{equation}\label{BB}
\lim_{n \to +\infty} \big\vert  \phi(Z_n(\omega))-\lambda^{-1}\phi(A_s\cdot Z_n(\omega))\big\vert^2 q(Z_n(\omega))=0,
\end{equation}
\begin{equation} \label{CC}
\lim_{n \to +\infty} \big\vert  \phi(Z_n(\omega))-\lambda^{-2}\phi(H_tH_s\cdot Z_n(\omega))\big\vert^2 p(Z_n(\omega))p(H_s\cdot Z_n(\omega))=0,
\end{equation}
and
\begin{equation}\label{DD}
\lim_{n \to +\infty} \big\vert
\phi(Z_n(\omega))-\lambda^{-2}\phi(A_tA_s\cdot
Z_n(\omega))\big\vert^2 q(Z_n(\omega))q(A_s\cdot Z_n(\omega))=0.
\end{equation}

 There are two cases to explore.
\begin{enumerate}
\item  \underline{  $\phi(0)=\phi(1)$}

Fix  $\omega \in \Omega_x$ and a cluster value $z_\omega$ of the sequence $(Z_n(\omega))_{n \geq 0}$.

If $p(z_\omega)\neq 0$, then, applying (\ref{AA}) with $ s $ arbitrarily close to $ 0 $, it yields
 $$
  \phi(z_\omega)=  \lambda^{-1}\phi(0).
 $$
 If  $p(z_\omega)= 0$,  we conclude  similarly with  (\ref{BB})    that  $ \phi(z_\omega)= \lambda^{-1}  \phi(1).$

 Consequently,  since $\phi(0)=\phi(1)$, the sequence
$(\phi(Z_n(\omega))_{n \geq 0}$  converges to $\Phi^\infty(\omega)
=  \lambda^{-1}\phi(0)$  and  $\phi(x) = \mathbb E_x(\Phi^\infty)=
\lambda^{-1}\phi(0)$. Thus, the function $\phi$ is    constant and
$\lambda =1$.

  \item  \underline{  $\phi(0)\neq \phi(1)$}

  Without loss of generality, we assume  $\phi(0)\neq 0$; the case $\phi(1)\neq 0$ is treated
   the same way.

 \begin{enumerate}
\item  First, assume that  there   exists $x \in [0, 1]$ and  $\omega_x\in \Omega_x$ such that the sequence  $(Z_n(\omega_x))_{n \geq 0}$ possesses a cluster point  $z_{\omega_x}$ with
 $p(z_{\omega_x})>0$.  Applying first (\ref{AA}) with $ s $ arbitrarily close to $ 0 $ and second (\ref{CC}) with  $s$  arbitrarily close to
 $1$ (so that  $p(H_s\cdot z_{\omega_x})>0$) and  $t$
 arbitrarily close to  $0$, it yields
  $$ \phi(z_{\omega_x})= \lambda^{-1}\phi(0) = \lambda^{-2}\phi(0). $$
 The condition   $\phi(0)\neq 0$ readily implies $\lambda =1$ and
thus $\phi \in \mathbb C \cdot {1\hspace*{-0.16cm}1}$.

 \item Assume  that $q(z_\omega)=1$  for any  $x \in [0, 1]$, any $\omega\in \Omega_x$ and any   cluster  values $z_\omega$ of the sequence $(Z_n(\omega))_{n\geq 0}$.

 Applying  first (\ref{BB})  with $s$ arbitrarily close to $0$ and second (\ref{DD}) with  $s$ arbitrarily close to  $1$ (so that $q(A_s\cdot z_{\omega})>0$) and   $t$ arbitrarily close to  $0$, it yields
   $$ \phi(z_{\omega})= \lambda^{-1}\phi(1) = \lambda^{-2}\phi(1). $$
If $\phi(1)\neq 0$, we deduce as above that  $\lambda =1$  and  $\phi \in \mathbb C \cdot  {1\hspace*{-0.16cm}1}$.
If   $\phi(1)=0$, the sequence   $(\lambda^{-n}\phi(Z_n(\omega))_{n \geq 0}$ converges to $0$ and the  martingale equality
$\phi(x)= \mathbb E_x(\lambda^{-n_l}\phi(Z_{n_l}))$ yields  $\phi \equiv 0$.
 \end{enumerate}
\end{enumerate}

Eventually, the operator $Q$ is quasi-compact on $\mathcal   H_\alpha[0,1]
$ with spectral radius equals $1$, its peripheral spectrum   is
reduced to $\{1\}$  and  the  characteristic subspace
associated to  $1$  equals ${\rm Ker}(Q-Id)$. More precisely, we  have the 4 following cases.

\begin{enumerate}
\item \underline{If  $ q(0)>0$ and $ p(1)>0$}, there exists  a bounded linear operator
    $R$ on $\mathcal   H_\alpha[0,1] $  with spectral radius  $\rho \in [0, 1[$ such that,
     for any  $\varphi \in \mathcal   H_\alpha[0,1] $ and   $n \geq 0$,
$$Q^n\varphi= \Bigl(\int_0^1\varphi(x)  f_p(x) {\rm d}x\Bigl) {1\hspace*{-0.16cm}1}+R^n\varphi.
$$
In this case, the chain $(Z_n)_{n\geq 0}$ is recurrent on $[0, 1]$.
\item \underline{If  $ q(0)=0$ and $ p(1)>0$}, there exists  a bounded operator  $R$ on $\mathcal   H_\alpha(E) $ with  spectral radius  $\rho \in [0, 1[$ such that, for any   $\varphi\in \mathcal   H_\alpha[0,1] $ and  $n \geq 0$,
\begin{equation} \label{delta0}
Q^n\varphi= \varphi(0) {1\hspace*{-0.16cm}1}+R^n\varphi.
\end{equation}
In this  case,  for any  $x \in [0, 1]$,  the chain $(Z_n)_{n\geq 0}$  converges  $\mathbb P_x$-a.s. to $0$; furthermore, for any $\epsilon \in ]0, 1[$,  the set $ [\epsilon, 1]$ is transient and there  exists $\kappa_\epsilon >0$ such that
$$ \mathbb P_x(Z_n\in  [\epsilon, 1])\leq \kappa_\epsilon \rho^n.
$$

\item\underline{  If  $ q(0)>0$ and $ p(1)=0$}, there exists  a bounded operator  $R$ on $\mathcal   H_\alpha[0,1] $ with spectral radius  $\rho \in [0, 1[$ such that, for any  $\varphi\in \mathcal   H_\alpha[0,1] $ and   $n \geq 0$,
\begin{equation*}
Q^n\varphi= \varphi(1) {1\hspace*{-0.16cm}1}+R^n\varphi.
\end{equation*}
For any   $x \in [0, 1]$,  the chain $(Z_n)_{n\geq 0}$  converges  $\mathbb P_x$-a.s. to  $1$; furthermore, for any $\epsilon \in ]0, 1[$,   the set $ [0, 1-\epsilon]$ is transient and there  exists $\kappa_\epsilon >0$ such that
$$ \mathbb P_x(Z_n\in  [0, 1-\epsilon ])\leq \kappa_\epsilon \rho^n.
$$

\item \underline{If $ q(0)=0$ and $ p(1)=0$},  there exists an harmonic  function   $h: [0, 1] \to [0, 1]$ such that $h(0)=0$ and $h(1)=1$  and   a bounded operator  $R$ on $\mathcal   H_\alpha[0,1] $ with spectral radius  $\rho \in [0, 1[$ such that, for any  $\varphi\in \mathcal   H_\alpha[0,1] $ and   $n \geq 0$,
\begin{equation*}
Q^n\varphi= \varphi(0) ({1\hspace*{-0.16cm}1}- h)+\varphi(1)  h+R^n\varphi.
\end{equation*}
For any   $x \in [0, 1]$,  the chain $(Z_n)_{n\geq 0}$  converges  to $0$ with probability $1-h(x)$ and to $1$ with  probability  $h(x)$. Indeed,   the bounded martingale $(h(Z_n))_{n \geq 0}$  converges $\mathbb P$-a.s. Since $h(0)\neq h(1)$,  it follows that    $(Z_n)_{ n\geq 0}$  converges  to a random variable  $Z_\infty$ with values in $\{0, 1\}$. The martingale property  yields $h(x)= \mathbb E_x(h(Z_\infty))=\mathbb P_x(Z_\infty = 1)$.
Consequently, for any $\epsilon \in ]0, 1[$,   the set $ [0, 1-\epsilon]$ is transient and there  exists $\kappa_\epsilon >0$ such that
$$ \mathbb P_x(Z_n\in  [0, 1-\epsilon ])\leq \kappa_\epsilon \rho^n.
$$
\end{enumerate}

\vspace{5mm}

    \noindent {\bf Example :  $  p(x)=1-x$.}

\noindent  We are in the   Case  4 above and the harmonic function
$h(x)=x$. In particular,   the sets  $[\epsilon, 1-\epsilon],
0<\epsilon <1$, are  transient.
\\
This transience property can  be obtained in a different way which is also of interest and  we
 present briefly.
Let us introduce the quantity $\Delta$ defined by: for any  $x$ in $[0, 1]$,
$$\Delta(x):= {\rm dist}(x, \{0, 1\})= \inf(x, 1-x).$$
Let us compute $\mathbb E_x(\Delta(Z_1))$.  We assume $x\in ]0, {1\over 2}]$,
the case $x\in [  {1\over 2}, 1[$ can be treated in a similar way.
\begin{eqnarray*}
\mathbb E_x(\Delta(Z_1))&=& {1-x\over x} \int_{0}^x  y\ {\rm d}y +
 {x\over1-x} \int_x^{1\over 2}y\  {\rm d}y+  {x\over1-x} \int_{1\over 2}^1(1-y) {\rm d}y
 \\
 &=& {3x-4x^2\over  4(1-x)}\\
 &\leq &{3\over 4} x.
 \end{eqnarray*}
Hence $\displaystyle \mathbb
E(\Delta(Z_n)|  \mathcal{F}_{n-1}) \leq {3\over 4}
\Delta(Z_{n-1})$ for any $n \geq 1$ and,  iterating,
$$
\forall n\geq 1, \forall x \in [0, 1]\qquad \mathbb E_x(\Delta(Z_n))\leq \left({3\over 4}\right)^n \Delta(x)\leq \left({3\over 4}\right)^n.
$$
Consequently  $\displaystyle \mathbb E_x\Bigl( \sum_{n=0}^{+\infty} \Delta(Z_n)\Bigr)=
\sum_{n=0}^{+\infty} \mathbb E_x\Bigl(  \Delta(Z_n)\Bigr)<+\infty$, so that  the sequence $(\Delta(Z_n))_{n\geq 1}$ converges $\mathbb P_x$-a.s. to $0$.


\normalsize

\begin{thebibliography}{99}

\bibitem{BDEG} Barnsley, M. F.,  Demko  S. G.,  Elton, J. H. \&  Geronimo, J. S.  (1988)  \emph{Invariant measures for Markov processes arising from iterated function systems with place-dependent probabilities}  Ann. Inst. H. Poincar\'e Probab. Statist.,  vol.  24, no. 3, 367--394.

\bibitem{BE} Barnsley, M. F. \& Elton, J.H.  (1988) \emph{A new class of Markov processes for image encoding}  Adv. in Appl. Probab.,  vol.   20, no. 1, 14--32.

\bibitem{BEH} Barnsley, M.F., Elton, J. H. \& Hardin, D. P. (1989) \emph{ Recurrent iterated function systems. Fractal approximation} Constr. Approx.,  vol.   5 , no. 1, 3--31.




\bibitem{DF}   Diaconis  P. \&  Freedman D.   (1999)
\emph{Iterated Random Functions}
SIAM review,  vol.  41, no. 1, 45--76.





\bibitem{DS}   Diaconis P. \&   Shahshahani M.  (1984)
\emph{ Products of random matrices as they arise in the study of
random walks on groups} Technical report no. 229, Departement of
Statistics, Stanford
University.




\bibitem{E} Elton, John H.  (1987)\emph{An ergodic theorem for iterated maps}  Ergodic Theory Dynam. Systems,  vol.  7 , no. 4, 481--488.




\bibitem{GH} Guivarc'h Y. \& Hardy J.   (1988)\emph{Th\'eor\`emes limites pour une classe de cha\^nes de Markov
et applications aux diff\'eomorphismes d'Anosov.} Ann. Inst. Henri
Poincar\'e,   vol.   24,  no. 1, p. 73--98.

\bibitem{GR} Guivarc'h Y. \& Raugi   (1986) \emph{Products of random matrices and
  convergence theorems,} Contempory Math., A.M.S.,  vol.   50,   31--54.




 \bibitem{H}  Hennion  H.  (1993) \emph{ Sur un th\'eor\`eme spectral et son application aux noyaux
Lipschitziens,} Proceedings of the American Mathematical Society,
  vol.  118, 627--634.

\bibitem{HH}  Hennion  H. \&   Herv\'e L.  (2001) \emph{ Limit theorems for
Markov chains and stochastic properties of dynamical systems by
quasi-compacity,} Lecture Notes in Math.,  vol.   1766, Springer.

\bibitem{Her} Herv\'e L.  (1994) \emph{ \'Etude d'op\'erateurs quasi-compacts positifs. Applications aux op\'erateurs de transfert,} Ann. Inst. Henri Poincar\'e,  section B,   vol.  30, no.
3,    437--466.


\bibitem{ITM} Ionescu-Tulcea  C.T. \&    Marinescu G.  (1950) \emph{Th\'eorie ergodique pour des classes d'op\'erations
non compl\`etement continues}, Annals of Mathematics,   vol.  52,
140--147.


\bibitem{L}   Letac G. (1986)
 \emph{A contraction principle for certain Markov Chains and its
 applications.}
  Contemporary Mathematics,    vol.  50  
    American Mathematical Society, Rhodes Island, N-Y.




\bibitem{MR}  Mirek M.  (2011) \emph{Heavy tail phenomenom and convergence to stable laws for iterated Lipschitz random maps}\\
Probability Theory and Related Fields,  vol.  151, 705--734.

\bibitem{P1} Peign\'e M. (1993) \emph {Iterated function systems
and spectral decomposition of the associated Markov operator}, in
Fascicule de probabilit\'es, Publ. Inst. Rech. Math. Rennes,  vol.
1993, Univ. Rennes I,  http://www.lmpt.univ-tours.fr/~peigne/fichiers/ifs.pdf.


\bibitem{P2}   Peign\'e  M.  (1992) \emph{Marches de Markov sur le semi-groupe des contractions de $\mathbb R^d$. Cas de la marche al\'eatoire \`a pas markoviens sur $(\mathbb R^+)^d$ avec chocs \'elastiques sur les axes},
Ann. Inst. Henri Poincar\'e, section B,   28, no.
1,  63--94.

 \bibitem{ST}  Stenflo \"O,    (2012) \emph{A survey of average contractive iterated function systems} J.
  Difference Equ. Appl. 18, no.
8, 1355--1380.

\bibitem{R}  Revuz  A.  (1975) \emph{Markov Chains}
North Holland.
\end{thebibliography}
\end{document}